\documentclass[12pt]{amsart}
\usepackage{amscd,amssymb}
\usepackage[graph,frame,poly,arc]{xy}  
\usepackage[plainpages,backref,urlcolor=blue]{hyperref}

\theoremstyle{plain}
\newtheorem{thm}[subsection]{Theorem}
\newtheorem{lem}[subsection]{Lemma}
\newtheorem{prop}[subsection]{Proposition}
\newtheorem{cor}[subsection]{Corollary}

\theoremstyle{definition}
\newtheorem{rk}[subsection]{Remark}
\newtheorem{definition}[subsection]{Definition}
\newtheorem{ex}[subsection]{Example}
\newtheorem{conj}[subsection]{Conjecture}

\numberwithin{equation}{section}
\newcommand{\Z}{\mathbb{Z}}
\newcommand{\Q}{\mathbb{Q}}

\newcommand{\C}{\mathbb{C}}
\newcommand{\PP}{\mathbb{P}}

\DeclareMathOperator{\rank}{rank}

\begin{document}

\title [Nearly free divisors and rational cuspidal curves]
{Nearly free divisors and rational cuspidal curves}

\author[Alexandru Dimca]{Alexandru Dimca$^1$}
\address{Univ. Nice Sophia Antipolis, CNRS,  LJAD, UMR 7351, 06100 Nice, France. }
\email{dimca@unice.fr}

\author[Gabriel Sticlaru]{Gabriel Sticlaru}
\address{Faculty of Mathematics and Informatics,
Ovidius University,
Bd. Mamaia 124, 900527 Constanta,
Romania}
\email{gabrielsticlaru@yahoo.com }
\thanks{$^1$ Partially supported by Institut Universitaire de France.}

\subjclass[2010]{Primary 14H45, 14B05; Secondary  14H50, 14C20}

\keywords{free divisor, rational cuspidal curve, Jacobian ideal, Milnor algebra}

\begin{abstract} We define a class of plane curves which are close to the free divisors in terms of the local cohomology of their Jacobian algebras and such that conjecturally any rational cuspidal curve $C$ is either free or belongs to this class. Using a recent result by Uli Walther we prove this conjecture when the degree of $C$ is either even or a prime power, and when the group of $C$ is abelian. 

\end{abstract}
 
\maketitle


\section{Introduction} \label{sec:intro}

Let $f$ be a homogeneous polynomial  in the polynomial ring $S=\C[x,y,z]$ and denote by $f_x,f_y,f_z$ the corresponding partial derivatives.
Let $C$ be the plane curve in $\PP^2$ defined by $f=0$ and assume that $C$ is reduced and not a union of lines passing through one point. We denote by $J_f$ the Jacobian ideal of $f$, i.e. the homogeneous ideal of $S$ spanned by the partial derivatives $f_x,f_y,f_z$ and let $M(f)=S/J_f$ be the corresponding graded ring, called the Jacobian (or Milnor) algebra of $f$. Let $I_f$ denote the saturation of the ideal $J_f$ with respect to the maximal ideal ${\bf m}=(x,y,z)$ in $S$.  Consider the graded $S-$submodule $AR(f) \subset S^{3}$ of {\it all relations} involving the derivatives of $f$, namely
$$\rho=(a,b,c) \in AR(f)_m$$
if and only if  $af_x+bf_y+cf_z=0$ and $a,b,c$ are in $S_m$. We set  $m(f)_k=\dim M(f)_k$ for any integer $k$. 

The curve  $C$ is a free divisor if $AR(f)$ is a free graded $S$-module of rank two. Freeness in the local analytic setting was introduced by K. Saito in \cite{KS2} and has attracted a lot of interest in the last decades, see for instance the long reference list in \cite{N}.
More precisely, a divisor $C:f=0$ of degree $d$ is {\it free with exponents} $(d_1,d_2)$ if the minimal resolution of the Milnor algebra $M(f)$ has the following form, 
\begin{equation} \label{r1}
0 \to \oplus_{i=1,2}S(-d_i-(d-1)) \to S^3(-d+1) \to S.
\end{equation} 
It is known that in this case $d_1+d_2=d-1$ and $\tau(C)=(d-1)^2-d_1d_2.$ Here $d_i$ are the degrees of the generating relations $r_i$, for $i=1,2$ for the free $S$-module $AR(f)$,
and $\tau(C)$ is the total Tjurina number of $C$, that is
$\tau(C)=\sum_{i=1,p} \tau(C,{x_i})$, the $x_i$'s being the singular points of $C$.
It is also known, see \cite{ST}, \cite{DS14},  that $C:f=0$ is a free divisor if and only if
\begin{equation} \label{r1b}
N(f)=I_f/J_f=H^0_ {\bf m}(M(f))=0.
\end{equation} 

In the recent paper \cite{DStFD} we have shown
that the freeness of the plane curve $C$ can be characterized in terms of the dimensions $m(f)_k$ and that a large number of the rational cuspidal curves as classified in \cite{FLMN}, \cite{F}, \cite{FZ}, \cite{SaTo}
give examples of irreducible free divisors.  The aim of this note is to define a class of  curves, called {\it nearly free divisors}  (and a subclass of it formed by the  {\it almost free divisors}), which are not free, but which are very close to the free divisors. This class contains conjecturally all the rational cuspidal curves. We have checked in this note that it contains indeed all the rational cuspidal curves listed in the papers mentioned above, at least the beginning part (or some special members) of each series of such curves. On the other hand, cuspidal curves which are not rational, as the Zariski sextic with 6 cusps, or rational curves which are not cuspidal, as the nodal cubic, do not enter in this class, see Examples \ref{ex0} and \ref{exdeg3+4}. These examples and those in \cite{DStFD} motivate the following conjecture.

\begin{conj}
\label{q2}
\noindent (i) Any rational cuspidal curve $C$ in the plane is either free or nearly free. 

\noindent (ii) An irreducible plane curve $C$ which is either free or nearly free is  rational.
\end{conj}

Using a recent result by U. Walther in \cite{Wa}, we prove the claim (i) in Conjecture \ref{q2}
for all curves of even degree, see Theorem \ref{thm2}, which is the our main result. The proof also implies that  this conjecture holds for a curve $C$ with an abelian fundamental group $\pi_1(\PP^2 \setminus C)$ or having a degree a prime power, see Corollary \ref{corPi}. 

 Moreover, any unicuspidal rational curve with a unique Puiseux pair is either free or nearly free, see Corollary \ref{corTop}, except the curves of odd degree in one case of the classification of such unicuspidal curves obtained in  \cite{FLMN} and recalled in Theorem \ref{thmClass} below, when our methods do not apply.
Indeed, we need a topological assumption on the cusps which is not fulfilled all the time when the degree is odd, see  Theorem \ref{thm2}.

As for  the claim (ii) in Conjecture \ref{q2}, note that non-irreducible nearly free curves, exactly as the free curves, may have irreducible components which are not rational, see Example \ref{nonlinearr}.
A related conjectural property of (nearly) free irreducible curves is the following.

\begin{conj}
\label{q3}
\noindent (i) Any free irreducible plane curve $C$ has only singularities with at most two branches. 

\noindent (ii) Any nearly free irreducible plane curve $C$ has only singularities with at most three branches. 
\end{conj}

We  prove the claim (ii) in Conjecture \ref{q2} and Conjecture \ref{q3} in a very limited number of cases in Theorem \ref{thmmain}.
Nanduri \cite{N} has constructed free rational curves with a unique singularity, and this singularity has either one or two branches, see Example \ref{exNanduri}. We construct an almost free rational curve having a singularity with 3 branches in Example \ref{ex3bran}.

To define the class of nearly free curves $C:f=0$, we allow a slightly more complicated minimal resolution for the 
Milnor algebra $M(f)$, compare \eqref{r2} below to \eqref{r1}, which follows from the condition $\dim N(f)_k \leq 1$ for all $k$, see Definition \ref{def} and Theorem \ref{thm0}.  The fact that the curve $C$ is almost free or nearly free depends only on the dimensions $m(f)_k$ of the homogeneous components of the Milnor algebra $M(f)$, see Corollary \ref{corC}.

As in the case of free curves, the Tjurina number $\tau(C)$ and other important numerical invariants of the nearly free curve $C:f=0$ are determined by the exponents $(d_1,d_2,d_3)$, which satisfy the properties $d_2=d_3$  (in view of this we call the pair $(d_1,d_2)$ the exponents of $C$) and  $d_1+d_2=d$, the degree of $C$, see Theorem \ref{prop1}. In Corollaries \ref{corA} and  \ref{corB} we present two possibilities for the exponents $(d_1,d_2)$ which cover many of  the encountered examples and are interesting since for them the existence of the resolution \eqref{r2} implies the nearly freeness. 

In the case of line arrangements, the point of view 
of our note might lead to progress on the main conjecture in the area, namely on whether the freeness of a (hyperplane) arrangement is determined by the combinatorics, see \cite{Te}, \cite{Yo}. We prove a property related to the factorization property of the characteristic polynomial of a free line arrangement in Proposition \ref{propD}, see also Example \ref{exchar} for a new look on a classical example in arrangement theory.

The computations of various invariants given in this paper were made using two computer algebra systems, namely CoCoA \cite{Co} and Singular \cite{Sing}, and they play a key role especially in the final section.
The corresponding codes are available on request, some of them being given in the book \cite{St}.

The authors thank Aldo Conca for a suggestion that lead us to simplify the proof of Theorem \ref{thm0}.

\section{The definition of nearly free and almost free divisors} 
We have seen that for a reduced curve  $C:f=0$  the existentce of a resolution \eqref{r1} is equivalent to the vanishing of the $S$-graded module $N(f)=I_f/J_f$, called the Jacobian module in \cite{Wa}, and the curves satisfying these equivalent properties are called free. The class of curves introduced in this note is defined by imposing the condition that the Jacobian module $N(f)$ is non-zero, but as small as possible.

\begin{definition}
\label{def}
A reduced plane curve $C:f=0$ is said to be {\it nearly free} if  $N(f) \ne 0$ and $\dim N(f)_k \leq 1$ for any $k$. We say that $C$ is {\it almost free} if $\dim N(f)=1$ for $d$ even, and $\dim N(f)=2$ for $d$ odd.
\end{definition}
Note that if $N(f) \ne 0$, then the values of  $\dim N(f)$ given for an almost free curve are the minimal possible ones, see Proposition \ref{prop2} below and its proof. The first result says that a nearly free divisor
has a minimal resolution of $M(f)$  slightly more complicated than that of a free divisor.

\begin{thm} \label{thm0}
If $C:f=0$ is a nearly free divisor, then  the minimal resolution of the Milnor algebra $M(f)$ has the form
\begin{equation} \label{r2}
0 \to S(-b-2(d-1)) \to \oplus_{i=1,3}S(-d_i-(d-1)) \to S^3(-d+1) \to S,
\end{equation} 
for some integers $1 \leq d_1\leq d_2 \leq d_3$ and $ b$.
 
\end{thm}

More information on this resolution is given in Theorem \ref{prop1} below.

\proof

By Hilbert Syzygy Theorem, see \cite{Eis}, p. 3, we know that the minimal resolution has the form
$$0 \to F_3 \to F_2 \to F_1=S^3(-d+1) \to F_0=S,$$
where $F_3= \oplus S(b_j)$ and $F_2= \oplus S(a_i)$ are graded free $S$ modules, with $a_i$, $b_j$ finite sets of integers. It is easy to see that $\rank F_2=\rank F_3+2$, since $M(f)$ has a constant Hilbert function, $m(f)_k=\tau(C)$ for $k$ large enough by  \cite{CD}.
Hence we have just to show that $\rank F_3=1$. 

If $F_i=\oplus S(-j)^{\beta _{i,j}}$, then the positive integers $\beta _{i,j}$ are called the graded Betti numbers of $M(f)$, they are denoted by $\beta _{i,j}(M(f))$ and they satisfy
$$\beta _{i,j}(M(f))=\dim Tor_i(M(f), \C)_j,$$
see \cite{Eis}, p. 8. Since $Tor(M(f), \C)=Tor(\C,M(f))$ and the graded $S$-module $\C$ has a minimal resolution given by the homological Koszul complex $C(x,y,z; S)$ of $x,y,z$ in $S$, it follows that
$$Tor_i(M(f), \C)=H_i(C(x,y,z; M(f)),$$
where $C(x,y,z; M(f))$ denotes the homological Koszul complex of $x,y,z$ in $M(f)$. 
The third differential in the complex $C(x,y,z; M(f))$ being given essentialy by multiplication by $x,y$ and $z$, it follows that $H_3(C(x,y,z; M(f)))$ is exactly the socle $s(M(f))$ of the module $M(f)$, that is 
$$s(M(f))=\{m \in M(f): \ \ xm=ym=zm=0 \}.$$
It follows that $\beta _{3,j}(M(f))=\dim s(M(f))_j$. On the other hand, it is clear that $s(M(f))=s(N(f))$, since we have the equality $N(f)=H^0_{\bf m}(M(f))$, the local cohomology of $M(f)$ with respect to the maximal ideal ${\bf m} =(x,y,z)$, see \cite{Se}. The module $N(f)$ is a module with a self dual resolution, see for instance \cite{DS1} or \cite{Se}.
If we set $\beta _{i}(M(f))= \sum_j \beta _{i,j}(M(f))$ and similarly for $N(f)$, this implies
$$\beta_3(M(f))=\beta_3(N(f))=\beta_0(N(f)).$$
In other words, $\beta_3(M(f))=1$, which is the claim of the theorem, is equivalent to $\beta_0(N(f))=1$, i.e. $N(f)$ is a cyclic $S$-module. But this last claim is a direct consequence of Corollary 4.3 in \cite{DPop}, which says that $N(f)$ has a Lefschetz type property and of our assumption that $\dim N(f)_j \leq 1$ for any $j$.

\endproof

\begin{rk}
\label{rkreso}

To construct  a resolution \eqref{r2} for a given polynomial $f$ we need the following ingredients.

\medskip

\noindent (i) Three relations $r_i=(a_i,b_i,c_i) \in S^3_{d_i}$ for $i=1,2,3$ among $f_x,f_y,f_z$, i.e. $$a_if_x+b_if_y+c_if_z=0,$$
necessary to construct the morphism 
$$\oplus_{i=1,3}S(-d_i-(d-1)) \to S^3(-d+1)$$ 
by the formula
$$(u_1,u_2,u_3) \mapsto u_1r_1+u_2r_2+u_3r_3.$$

\noindent (ii) One relation $R=(v_1,v_2,v_3) \in \oplus_{i=1,3}S(-d_i-(d-1))_{b+2(d-1)} $ among $r_1,r_2,r_3$, i.e. $v_1r_1+v_2r_2+v_3r_3=0$, necessary to construct the morphism 
$$S(-b-2(d-1)) \to \oplus_{i=1,3}S(-d_i-(d-1))$$
by the formula $w \mapsto wR$. Note that $v_i \in S_{b-d_i+d-1}$.

\end{rk}

\begin{ex} \label{ex1}
If $C$ is a smooth curve, then we have a resolution as in \eqref{r2} with $d_1=d_2=d_3=d-1$. Hence $b=d-1$ and $\tau(C)=0$. This example explains the normalization choosen for $b$.
In this case the resolution \eqref{r2} is just the Koszul complex of the partial derivatives $f_x,f_y,f_z$. We can take $r_1=(f_y,-f_x,0)$, $r_2=(f_z,0,-f_x)$, $r_3=(0,f_z,-f_y)$ and $R=(f_z,-f_y,f_x)$. Note that in this case $I_f/J_f=M(f)$ is far from being trivial as in \eqref{r1b}.
Moreover $\tau(C)=0$ is in sharp contrast to the property 
$$\tau(C) \geq \frac{3}{4}(d-1)^2$$
valid for an irreducible  free divisor $C$, see \cite{DStFD}, Theorem 2.8. In conclusion the class of curves defined by having a resolution of type \eqref{r2} is too large.
\end{ex}

\begin{ex} \label{ex0}
Consider the famous Zariski sextic curve with six cusps on a conic $C:f=(x^2+y^2)^3+(y^3+z^3)^2$. Then the minimal resolution of the Milnor algebra $M(f)$ has the form
$$0 \to  S(-11) \oplus  S(-12) \to  S(-8) \oplus S(-10)^3 \to  S(-5)^3 \to S,$$
and hence it is not of the type described in \eqref{r2}. 

The line arrangement $C: f=(y-z)(y-2z)(y-x)(y+x)(x+y+z)$ consists of 5 lines in general position.
It has the minimal resolution of the Milnor algebra $M(f)$ of the form
$$0 \to  S(-8)^2 \to  S(-7)^4 \to  S(-4)^3 \to S,$$
hence again not of the type described in \eqref{r2}. 

It follows that, by imposing a minimal resolution of 
type \eqref{r2} for the Milnor algebra $M(f)$, we get indeed a strict subclass in the class of reduced plane curves, and even in the class of line arrangements.
\end{ex}

The first natural question is whether such nearly free curves exist. The following examples show that the answer is positive.

\begin{ex} \label{exdeg3+4}
In degree $d=3$, consider a conic plus a secant line, e.g. $C: x^3+xyz=0$. Then $C$ is almost free with the resolution for $M(f)$ of the form
$$0 \to S(-5) \to S(-3) \oplus S(-4)^2 \to S(-2)^3 \to S,$$
and hence $(d_1,d_2,d_3)=(1,2,2)$. On the other hand, the nodal cubic $C:f=xyz+x^3+y^3=0$
is not nearly free, since $\dim N(f)_1=2$.

In degree $d=4$ we consider first the quartic with 3 cusps, e.g. 
$$C: x^2y^2+y^2z^2+x^2z^2-2xyz(x+y+z)=0.$$
Then $C$ is almost free with the minimal resolution for $M(f)$ of the form
$$0 \to S(-6) \to  S(-5)^3 \to S(-3)^3 \to S,$$
and hence $(d_1,d_2,d_3)=(2,2,2)$.
The quartic $C: z^4-xz^3-2xyz^2+x^2y^2=0$ (resp. $C: y^4-2xy^2z+yz^3+x^2z^2=0$) has two cusps of type $A_2$ and $A_4$ (resp. one cusp of type $A_6$), and it is almost free with the same resolution for $M(f)$ as for the 3-cuspidal quartic.

\end{ex} 

\begin{ex} \label{linearr} We give now some examples of line arrangements which are almost free.

\noindent (i) Let $C: f=xyz(x+y+z)=0$ be the union of 4 lines in general position. Then $C$ is an almost free curve with the resolution for $M(f)$ given by
$$0 \to S(-6) \to S(-5)^3 \to S(-3)^3 \to S,$$
with $(d_1,d_2,d_3)=(2,2,2)$ and $\dim N(f)_3=1$, and  $N(f)_k=0$ for other $k$'s.

\noindent (ii) Let $C: f=xyz(x-z)(x+y+z)=0$ be the union of 5 lines with a triple point $(0:1:0)$ and 7 nodes. Then we get an almost free curve with the resolution for $M(f)$ given by
$$0 \to S(-8) \to S(-6)\oplus S(-7)^2\to S(-4)^3 \to S$$
and hence $(d_1,d_2,d_3)=(2,3,3)$ with $\dim N(f)_k=1$ for $k=4,5$  and  $N(f)_k=0$ for other $k$'s.

\noindent (iii) Let $C: f=xyz(x-z)(y+2z)(x+y+z)=0$ be the union of 6 lines, with 3 triple points and 6 double points. Then we get an almost free curve with the resolution for $M(f)$ given by
$$0 \to S(-9) \to S(-8)^3 \to S(-5)^3 \to S$$
and hence $(d_1,d_2,d_3)=(3,3,3)$ with $\dim N(f)_6=1$, and  $N(f)_k=0$ for other $k$'s.

\end{ex} 

\begin{ex} \label{nonlinearr} 

(i) We give now an example of a non linear arrangement which is almost free and with one component which is  rational, but not cuspidal.
The curve $C: f=z((x^2+y^2+z^2)^3-27x^2y^2z^2)=0$ is a union of a line $L:z=0$ and a sextic curve $C'$ with 6 cusps (located at points of type $(0:1:\pm 1)$ and 4 nodes (located at the points $(1:\pm 1: \pm 1)$) as singularities. Hence $\delta (C)=10$, showing that $C'$ is  rational. The curve $C$ is an almost free curve with the resolution for $M(f)$ given by
$$0 \to S(-11) \to S(-9)\oplus S(-10)^2\to S(-6)^3 \to S$$
and hence $(d_1,d_2,d_3)=(3,4,4)$ with $\dim N(f)_k=1$ for $k=7,8$  and  $N(f)_k=0$ for other $k$'s.

(ii) Here is an example of a non linear arrangement which is nearly free and with one component which is a smooth cubic, hence not rational.

The sextic $C: f=(x^3+y^3+z^3)(x+y)(x+z)(y+z)=0$ is a union of 3 lines and an elliptic curve. The curve $C$ is a nearly free curve with the resolution for $M(f)$ given by
$$0 \to S(-10) \to S(-7)\oplus S(-9)^2\to S(-5)^3 \to S$$
and hence $(d_1,d_2,d_3)=(2,4,4)$, $b=0$,  with $\dim N(f)_k=1$ for $k=5,6,7$  and  $N(f)_k=0$ for other $k$'s.

\end{ex}

\section{Properties of  nearly free and almost free divisors}

We recall first the definition of some invariants associated with a Milnor algebra $M(f)$, see \cite{DStEdin}.

\begin{definition}
\label{def0}

For a reduced plane curve $C:f=0$ of degree $d$, three integers are defined as follows.

\noindent (i) the {\it coincidence threshold} 
$$ct(f)=\max \{q:\dim M(f)_k=\dim M(f_s)_k \text{ for all } k \leq q\},$$
with $f_s$  a homogeneous polynomial in $S$ of degree $d$ such that $C_s:f_s=0$ is a smooth curve in $\PP^2$.

\noindent (ii) the {\it stability threshold} 
$st(f)=\min \{q~~:~~\dim M(f)_k=\tau(C) \text{ for all } k \geq q\}.$

\noindent (iii) the {\it minimal degree of a  nontrivial (or essential) syzygy} 
$$mdr(f)=\min \{q~~:~~ H^2(K^*(f))_{q+2}\ne 0\},$$
where $K^*(f)$ is the Koszul complex of $f_x,f_y,f_z$ with the natural grading.
\end{definition}
Note that one has  for $j<d-1$ the following equality
\begin{equation} 
\label{ar=er}
AR(f)_j=H^2(K^*(f))_{j+2}.
\end{equation} 

It is known that one has
\begin{equation} 
\label{REL}
ct(f)=mdr(f)+d-2.
\end{equation} 
 Let $T=3(d-2)$ denote the degree of the socle of the Gorenstein ring  $M(f_s)$.
Now we  prove some results on nearly free curves.
\begin{thm} \label{prop1}
Suppose the curve $C:f=0$ has a minimal resolution for $M(f)$ as in \eqref{r2} with $d_1 \leq d_2 \leq d_3$. Then one has the following.

\noindent (i)  
$d_1+d_2 \geq d$,  $b=\sum_{i=1,3}d_i-2(d-1)$ and 
$$\tau(C)=(d-1)\sum_{i=1,3}d_i- \sum_{i < j}d_id_j.$$
Moreover  $mdr(f)=d_1$, $ct(f)=d_1+d-2$ and $st(f) =b+2d-4$.

\noindent (ii)  Suppose in addition that the curve $C:f=0$ is nearly free. Then one has the following:
$d_1+d_2 =d$, $d_2=d_3$,  $b=d_2-d+2$ and 
$$\tau(C)=(d-1)^2-d_1(d_2-1)-1.$$
Moreover $st(f) =d_2+d-2$ and $ct(f)+st(f)=T+2=3d-4$. The curve $C$ is almost free exactly when
$d_2=d_1$ for $d$ even, or $d_2=d_1+1$ for $d$ odd.
\end{thm}

In view of this result, we call the pair $(d_1,d_2)$ the {\it exponents of the nearly free curve} $C$.

\proof First we prove the part (i). The claims for $mdr(f)$ and $ct(f)$ are obviuos. The relation $r_2$ is not a multiple of $r_1$ (otherwise the resolution \eqref{r2} is not minimal). Then the claim $d_1+d_2 \geq d$ follows from Lemma 1.1 in \cite{ST}.

For any integer $j$, the resolution \eqref{r2} yields an exact sequence
\begin{equation} \label{r2j}
0 \to S_{j-b-2(d-1)} \to \oplus_{i=1,3}S_{j-d_i-(d-1)} \to S^3_{j-d+1} \to S_j \to M(f)_j \to 0.
\end{equation} 
For large $j$, this gives the equality
$\tau(C)= P(d,d_1,d_2,d_3,b;j)$, where the polynomial   $P(j)=P(d,d_1,d_2,d_3,b;j) \in \Q[j]$ is defined by
$$P(j)= {j+2 \choose 2} -3{j-d+3 \choose 2}+\sum_{i=1,2,3}{j-d_i-d+3 \choose 2}-{j-b-2d+4 \choose 2}.$$
The coefficient of $j^2$ in the right hand side is zero, the vanishing of the coefficient of $j$ gives the formula for $b$, while the constant term gives the value for $\tau(C)$, since $m(f)_j=\tau(C)$ for large $j$ by \cite{CD}.
For the claim on $st(f)$, note that for $j \geq b+2d-4$ the above binomial coefficients give the dimensions of the corresponding spaces in the exact sequence \eqref{r2j}, but not for $j=b+2d-5$, when the last term is 1 instead of being $0=\dim S_{-1}$, due to the equality
$${-1 \choose 2}=(-1)(-2)/2=1.$$

Now we prove the part (ii).
Consider the graded dual of the resolution \eqref{r2} twisted by $(-3)$, namely 
$$S \to S^3(d-4) \to  F_2'=\oplus S(-a_i-3) \to F_3'= S(b+2d-5) \to 0.$$
 Then the cokernel $Q$ of the morphism
$$ \delta: \oplus_iS(d_i+d-4) \to S(b+2d-5)$$
is the graded dual of $N(f)=H^0_{\bf m}(M(f))$ by Theorem A 1.9 in \cite{Eis}, p. 193. 
Indeed, $N(f)=H^0_{\bf m}(M(f))$, the local cohomology with respect to the maximal ideal ${\bf m} =(x,y,z)$, see \cite{Se}.
It follows that
$$Q=S/(v_1,v_2,v_3)(b+2d-5),$$
where $v_1,v_2,v_3$ are the homogeneous polynomials from Remark \ref{rkreso}.
This quotient is finitely dimensional if and only if $v_1,v_2,v_3$ is a regular sequence, and then the Hilbert function of the quotient depends only on the degree of the $v_i$'s. In particular, we can take $v_1=x^p$, $v_2=y^q$ and $v_3=z^r$ for some integers $p,q,r >0$ (when the quotient is the tensor product of the very simple rings $\C[x]/(x^p)$,  $\C[y]/(y^q)$, $\C[z]/(z^r)$.
In particular, to get at most one dimensional homogeneous components we need to have (up to a permutation) $\deg v_2=\deg v_3=1$ and $\deg v_1=k$ a positive integer.
It follows from Remark \ref{rkreso} that $b-d_i+d-1=1$ for $i=2,3$, and hence $d_2=d_3=b+d-2=d_1+2d_2-d$,
which implies $d_1+d_2=d$. The equality $b-d_1+d-1=k$ yields the last claim.

\endproof

\begin{cor}
\label{corA}
Suppose the curve $C:f=0$ of even degree $d=2e$ has a minimal resolution for $M(f)$ as in \eqref{r2} with $d_1 = d_2 =d_3=e$. Then the curve $C$ is almost free, $b=2-e$, $\tau(C)=3e(e-1)$, $mdr(f)=e$ and $ct(f)=st(f)=3e-2$.
\end{cor}

 \proof
The exact sequence \eqref{r2j} becomes in this case the following sequence
$$0 \to S_{j-3e} \to S_{j-3e+1} ^3\to S^3_{j-2e+1} \to S_j \to M(f)_j \to 0.$$
The associated polynomial $P(j)=P(d,d_1,d_2,d_3,b;j)$ is in this case
$$P(j)= {j+2 \choose 2} -3{j-2e+3 \choose 2}+3{j-3e+3 \choose 2}-{j-3e+2 \choose 2}.$$
We know that $P(j)= \tau(C)$ for any $j$. On the other hand, by writing the same exact sequence for a smooth curve $C_s:f_s=0$ of degree $d=2$, we see that 
$m(f)_k=m(f_s)_k$ for any $k\leq 3e-2=T/2+1$ where $T=3d-6$.

Note also that for $j \geq 3e-2$, the above binomial coefficients give the dimensions of the corresponding spaces in the exact sequence. In particular $m(f)_j= \tau(C)$ for $j \geq 3e-2$.
When we set $j=e-3=T/2$ we get $\tau(C)=P(3e-3)=m(f)_{3e-3}-1$ exactly as at the end of the proof of the claim (i) in Theorem \ref{prop1} above.

It remains now to use the formula
\begin{equation} \label{n(f)}
\dim N(f)_k= m(f)_k+m(f)_{T-k}-m(f_s)_k-\tau(C),
\end{equation} 
which is part of Corollary 3 in \cite{DS1}, but can be also obtained as follows.
Note that $\dim N(f)_k=m(f)_k -\dim (S/I_f)_k=m(f)_k-(\tau(C)-def_k\Sigma_f)$, with $def_k\Sigma_f$ the defect of the singular subscheme $\Sigma_f$ with respect to degree $k$ polynomials,
and $def_k\Sigma_f= m(f)_{T-k}-m(f_s)_k$, see for all this \cite{DBull}.
\endproof

The same type of argument with some tedious computations implies the following.

\begin{cor}
\label{corB}
Suppose the curve $C:f=0$ of  degree $d$ has a minimal resolution for $M(f)$ as in \eqref{r2} with $d_1 =1$ and $ d_2 =d_3=d-1$. Then the curve $C$ is nearly free,  more precisely $\dim N(f)_k=1$ for $d-2 \leq k \leq 2d-4$ and $N(f)_k=0$ otherwise. In addition one has $b=1$, $\tau(C)=(d-1)(d-2)$, $mdr(f)=1$, $ct(f)=d-1$ and $st(f)= 3d-3$.

\end{cor}

The formula \eqref{n(f)}  and Theorem \ref{prop1} also yield the following.

\begin{cor}
\label{corC}
The fact that a reduced curve $C$ is almost free or nearly free depends only on the dimensions $m(f)_k$ of the homogeneous components of the Milnor algebra $M(f)$. Conversely, two nearly free curves $C_1:f_1=0$ and $C_2:f_2=0$ with the same degree and the same total Tjurina numbers have the same exponents. In particular they satisfy $m(f_1)_k=m(f_2)_k$ for any $k$.
\end{cor}

Note that the other invariants of two such curves $C_1$ and $C_2$ can be quite different, see Propositions \ref{ST1} and \ref{ST2} below.

\begin{rk}
\label{rkexample}

In most of the examples we have encountered so far, the minimal resolution for $M(f)$ is of the type described in Corollary  \ref{corA} or in Corollary \ref{corB}.
Indeed, Corollary \ref{corA} applies to the three rational quartics in Example \ref{exdeg3+4}, the line arrangements (i) and (iii) in Example \ref{linearr}, to all almost free curves in Example \ref{ex2}, to the curves in Example \ref{prop0} and to the curves in the formula \eqref{cjk},
for $2 \leq k+j \leq 10$ from  Example \ref{ex3}.
Corollary \ref{corB} applies to the cubic in Example \ref{exdeg3+4}, to all nearly free curves in 
Proposition \ref{ST1} and Proposition \ref{ST2}.

\end{rk}

The formula \eqref{n(f)}  and Theorem \ref{prop1} also give the following related result.
Recall just that projective rigidity is equivalent to $N(f)_d=0$, see \cite{Se}.
\begin{cor}
\label{corC1}
Let $C:f=0$ be a nearly free curve of degree $d$ with exponents $(d_1,d_2)$. Then  $N(f)_k \ne 0$ for $ d+d_1-3 \leq k \leq d+d_2-3$ and $N(f)_k=0$ otherwise.
The curve $C$ is projectively rigid if and only if $d_1 \geq 4$. In particular, for $d\geq 8$, an almost free curve is rigid.
\end{cor}

\begin{prop} \label{prop2} Consider the reduced plane curve $C:f=0$ of degree $d$, which is not free.

\noindent (i)  If $d$ is even, then  $\dim N(f)_{T/2}=1$ if and only if $\dim N(f)_k \leq 1$ for all $k$.
Moreover, $\dim N(f)=1$ if and only if  $\dim N(f)_{T/2}=1$ and $\dim N(f)_{T/2-1}=0$. 

\noindent (ii) If $d$ is odd, then  $\dim N(f)_{[T/2]}=1$ if and only if $\dim N(f)_k \leq 1$ for all $k$.
Moreover, $\dim N(f) \geq 2$  and  the minimal value $\dim N(f)=2$ is attained  if and only if  $\dim N(f)_{[T/2]}=1$ and $\dim N(f)_{[T/2]-1}=0$.

\end{prop} 

\proof
If $N(f) \ne 0$, then there is a nonnegative integer $n(C)$,
$0 \leq n(C) \leq [T/2]$, such that $ N(f)_k=0$ for $k<n(C)$ or $k >T-n(C)$, and the remaining dimensions $\dim N(f)_k$ form a unimodal sequence of strictly positive numbers, see Corollary 4.3 in \cite{DPop},  symmetric with respect to the middle point $[T/2]$, see \cite{DS1} and \cite{Se}.
 This implies both claims above. 

\endproof

The next result and the following two examples look at the number of branches a singularity of a (nearly) free curve can have.
Let $\nu: \tilde C \to C$ the normalisation of the irreducible curve $C$. Assume that the singularities of $C$ are the points $p_i$, $i=1,...,m$ and the germ $(C,p_i)$ has $r_i$ branches.
\begin{thm}
\label{thmmain}
\noindent (i) If $C$ is an irreducible free curve of degree $d $, then 
$$2g(\tilde C)+\sum_{i=1,m}(r_i-1) \leq \frac{(d-1)(d-5)}{4}.$$
In particular, if $d \leq 6$, then $C$ is rational and $r_i\leq 2$ for all $i$, with equality for at most one $i$.

\noindent (ii) If $C$ is an irreducible nearly free curve of degree $d$ such that $d_1=1$ (resp. $d \leq 5$), then  $C$ is rational cuspidal (resp. rational and $r_i\leq 2$ for all $i$, with equality for at most one $i$.).

\end{thm}

\proof  Using the usual additivity of the Euler number $E$, we get
$$E(C)=E(C_{smooth})+m=E(\tilde C \setminus \nu^{-1}(C_{sing}))+m=2-2g(\tilde C)-\sum_{i=1,m}(r_i-1),$$
since $\nu^{-1}(p_i)$ consists of $r_i$ points. On the other hand one has
$$E(C)=2-(d-1)(d-2)+\mu(C),$$
where $\mu(C)$ is the sum of all local Milnor numbers $\mu(C,p_i)$, see for instance \cite{D1}.
It follows that for any irreducible curve $C$ one has
\begin{equation} \label{muC}
\mu(C)=(d-1)(d-2)-2g(\tilde C)-\sum_{i=1,m}(r_i-1).
\end{equation} 
This gives
\begin{equation} \label{gC}
2g(\tilde C)+\sum_{i=1,m}(r_i-1) =(d-1)(d-2)-\mu(C) \leq (d-1)(d-2)-\tau(C),
\end{equation} 
since $\tau(C) \leq \mu(C)$.
Suppose now first that $C$ is free. Then we have show in \cite{DStFD} that one has
$$\tau(C) \geq \frac{3}{4}(d-1)^2.$$
This implies
$$2g(\tilde C)+\sum_{i=1,m}(r_i-1) \leq (d-1)(d-2)-\frac{3}{4}(d-1)^2=\frac{(d-1)(d-5)}{4},$$
which yields the first claim. Assume now that $C$ is nearly free. The formula for $\tau(C)$ given 
in Proposition \ref{prop1} and the inequality \eqref{gC} yield
$$2g(\tilde C)+\sum_{i=1,m}(r_i-1) \leq (d-1)(d-2) -(d-1)^2+(d_1-1)(d_2-1)+d_2=(d_1-1)(d_2-2).$$
When $d_1=1$ (resp. $d\leq 5$), the upper bound is zero (resp. is at most one), so the claim follows.

\endproof

\begin{ex} \label{exNanduri} We present here the examples of free divisors given in \cite{N}.
Let $d\geq 5$ and choose two nonnegative  integers $j,k$ such that $0 \leq j+k \leq [(d+1)/2]-3$.
Let $F_1(x,y)$ (resp. $F_2(x,y)$) be a homogeneous polynomial in $x,y$ of degree $j$ (resp. $d-j_1$ with $j_1=[d/2]+j+1$), such that $F_1$ is square free and both $F_1$ and $F_2$ contain the monomials $x^j$ and $y^j$
(resp. $x^{d-j_1}$ and $y^{d-j_1}$) with nonzero coefficients. Then the curve
$$C: f=x^{d-j}F_(x,y)+y^{j_1}F_2(x,y)+x^ky^{d-k-1}z=0$$
is a free irreducible divisor with exponents $(d_1,d_2)=(v,v)$ if $d=2v+1$ and $(d_1,d_2)=(v-1,v)$ if $d=2v$. This gives the Tjurina number $\tau(C)=3v^2$ in the case $d=2v+1$ and 
$\tau(C)=3v^2-3v+1$ in the case $d=2v$ as noticed in \cite{N}. The curve $C$ has a unique singularity located at $p=(0:0:1)$ and in the local coordinates $x/z$ and $y/z$ it is a Newton non-degenerate and commode singularity, hence one can easily compute the Milnor number
$$\mu(C)=2A-2d+1=d^2-3d+1,$$
where $A$ is the area under the Newton polygon, see \cite{AGV}.
$C$ is a rational curve as $z$ occurs only with exponent one, hence the formula \eqref{muC} implies
that the number of branches of $(C,p)$ is  always $r=2$.
\end{ex} 

\begin{ex} \label{ex3bran} Using a similar formula as for Nanduri's polynomials, we have found that the curve
$$C: f=x^{11}+y^{11}+x^2y^6(x+5y)^2z=0$$
is rational, irreducible and has just one singularity, again  at $p=(0:0:1)$. This curve is an almost free divisor with exponents $(d_1,d_2)=(5,6)$, $\tau(C)=74$, $\mu(C)=88$ and the formula \eqref{muC} implies
that the number of branches of $(C,p)$ is  $r=3$ in this case. Note that in this case the singularity $(C,p)$ is no longer Newton nondegenerate, so the Milnor number is computed using Singular.
A nearly free divisor with similar properties is given by
$$C: f=x^{13}+y^{13}+x^2y^8(x+5y)^2z=0.$$
Here  $(d_1,d_2)=(5,8)$, $\tau(C)=108$, $\mu(C)=130$ and the formula \eqref{muC} implies
that the number of branches of $(C,p)$ is  again $r=3$.

\end{ex} 

We end this section with a remark on nearly free line arrangements.
It is known that for a free arrangement, the characteristic polynomial has a nice factorization in terms of the exponents, see \cite{Yo} for an excellent survey. For a nearly free line arrangement we have the following  similar result.

\begin{prop}
\label{propD} Let $C$ be a nearly free  arrangement of $d$ lines in $\PP^2$ with exponents $(d_1,d_2)$ and denote by $U$ its complement.
Then the characteristic polynomial $$\chi (U)(t)=t^2 -b_1(U)t+b_2(U)$$ is given by
$$\chi (U)(t) -1=(t-d_1)(t-(d_2-1)).$$
\end{prop}

\proof For any line arrangement one has $b_1(U)=d-1$, hence in our case $b_1(U)=d_1+d_2-1.$
Note that $E(U)=E(\PP^2)-E(C)=3+d^2-3d-\mu(C)$. It follows that
$$b_2(U)=E(U)+d-2=(d-1)^2-\mu(C).$$
Since all the singularities of $C$ are homogeneous, it follows that $\mu(C)=\tau(C)$ and using the formula for $\tau(C)$ given in Theorem \ref{prop1} we get $b_2(U)=d_1(d_2-1)+1$

\endproof

\begin{ex} \label{exchar} (i) For the line arrangement in Example \ref{linearr} (i) we have
$d_1=d_2=2$ and hence $\chi(U)(t)-1=(t-1)(t-2)$.

\noindent (ii) Consider the line arrangement $C: f=xy(x-y)(x-2y)(x-3y)z(x+5y+7z)=0$,
which occurs essentially in Example 4.139 in \cite{OT} as an illustration of the fact that 
$\chi(U)(t)$ can factor even for non free arrangements. It turns out that this curve is a nearly free arrangement with exponents $(d_1,d_2)=(2,5)$ and hence
$$\chi(U)(t)=(t-2)(t-4)+1=(t-3)^3.$$

\end{ex}

\section{Local and global Milnor fiber monodromies}

We state now the main result of our paper.

\begin{thm}
\label{thm2}
Let $C:f=0$ be a rational cuspidal curve of degree $d$. Assume that either 

\begin{enumerate}

\item $d$ is even, or

\item $d$ is odd and for any singularity $x$ of $C$, the order of any eigenvalue $\lambda_x$ of the local monodromy operator $h_x$ is not $d$.

\end{enumerate}

Then $C$ is either a free or a nearly free curve.

\end{thm}

For examples of rational cuspidal curves not satisfying the assumption in this theorem, see the case (3) with $d$ odd in Theorem \ref{thmClass} below.

\proof

Theorem 4.3 and the explicit description of the injection stated there given in the beginning of section 4.2 by U. Walther in \cite{Wa}, yield an injection
$$N(f)_{2d-2-j} \to H^2(F, \C)_{\lambda},$$
with $j=1,2,...,d$, where $F: f(x,y,z)-1=0$ is the Milnor fiber associated to $C$ and the subscript $\lambda$ indicates the eigenspace of the monodromy action corresponding to the eigenvalue $\lambda= \exp(2\pi i (d+1-j)/d)$. 

Suppose now that $d$ is even, say $d=2d_1$.
In view of Proposition \ref{prop2}, it is enough to show that
$\dim H^2(F,\C)_{\lambda} =1$ for  ${\lambda}=-1$ which corresponds to $j=d_1+1$.

Note first that $-1$ is not an eigenvalue for the local monodromy $h_x$ of any cusp $(C,x)$. Indeed, if this cusp has as Puiseux pairs the sequence
$(m_1,n_1)$,...,$(m_g,n_g)$ for some $g \geq 1$, it was shown by L\^e in \cite{Le} that the characteristic polynomial of the monodromy operator $h_x$ is given by 
\begin{equation} \label{AP}
\Delta(t)=P_{\ell_1,n_1}(t^{\nu_2}) \cdots P_{\ell_g,n_g}(t^{\nu_{g+1}}),
\end{equation}  
where  $\ell_1=m_1$,
$\ell_i=m_i+n_i(\ell_{i-1}n_{i-1}-m_{i-1})$ for $i=2,...,g$, and $\nu _i=n_i\cdots n_g$ with $\nu_{g+1}=1$ and 
$$P_{\ell,n}(t)=\frac{(t^{\ell n}-1)(t-1)}{(t^{\ell}-1)(t^{n}-1)}.$$
Using this formula, it is clear that $P_{\ell,n}(\pm 1) \ne 0$, and hence any factor $P_{\ell_j,n_j}(t^{\nu_{j+1}})$ does not admit $-1$ as a root.

 The fact that $\lambda=-1$ is not an eigenvalue for the local monodromy $h_x$ of any cusp $(C,x)$ implies that $H^1(F)_{\lambda}=0$, see Theorem 6.4.13 in \cite{D2}. 
Since $E(U)=E(\PP^2)-E(C)=1$, it follows that 
$$\dim H^2(F,\C)_{\lambda} -\dim H^1(F,\C)_{\lambda} +\dim H^0(F,\C)_{\lambda}=1.$$
Since clearly  $H^0(F)_{\lambda}=0$, the result is proved.

Suppose now that $d$ is odd, say $d=2d_1+1$. By Proposition \ref{prop2} again, it is enough to show that
$\dim H^2(F,\C)_{\lambda} =1$ for  ${\lambda= \exp(2\pi i d_1/(2d_1+1))   }$ which corresponds to $j=d_1+2$. As this eigenvalue ${\lambda}$ has clearly order $d$, we conclude as in the previous case.

\endproof

\begin{cor}
\label{corPi}
Let $C:f=0$ be a rational cuspidal curve  of degree $d$ such that

\begin{enumerate}

\item either $d=p^k$ is a prime power, or

\item $\pi_1(U)$ is abelian, where $U=\PP^2 \setminus C$.

\end{enumerate}

Then $C$ is either a free or a nearly free curve.
\end{cor}

\proof Using Theorem 1.1 in \cite{DL}, we get with obvious notation
$$\dim H^1(F, \C)_{\ne 1} =\dim H^3(X_C,\C),$$
where $X_C$ is the surface in $\PP^3$ given by $f(x,y,z)-t^d=0$.
Let $\rho: Y_C\to X_C$ be  a minimal resolution of singularities for $X_C$. It is known that $\rho^*: H^3(X_C,\C) \to H^3(Y_C,\C)$ is an isomorphism, see \cite{BK}, p.122.
Hence, using the proof of   Theorem \ref{thm2}, it is enough to show that $b_3(Y_C)=0$.
By Poincar\'e duality we know that $b_3(Y_C)=b_1(Y_C)$, and hence it is enough to show that $b_1(Y_C)=0$. When $d$ is a prime power, this fact goes back to Zariski \cite{Za}.

Note next that the Milnor fiber $F$ is a Galois degree $d$ cover of $U$, and if 
$\pi_1(U)$ is abelian, which actually means $\pi_1(U)=H_1(U)=\Z/d\Z$, it follows that $\pi_1(F)=0$. Since $F$ is an open subset in the smooth part of the surface $X_C$, it follows that $F$ can be regarded as a Zariski open subset of $Y_C$. This clearly implies $\pi_1(Y_C)=0$, thus completing the proof.

\endproof

Note that a lot of examples of rational cuspidal curves $C$ with an abelian fundamental group $\pi_1(U)$ can be found in the papers \cite{Art}, \cite{Ulu}. Here is one example of an infinite family of such curves.

\begin{ex}
\label{exPi}
Let $C$ be a rational cuspidal curve of type $(d,d-2)$ having three cusps. Then there exists a unique pair of integers $a,b$, $a\geq b \geq 1$ with $a+b=d-2$ such that up to projective equivalence the equation of $C$ can be written in affine coordinates $(x,y)$ as
$$f(x,y)= \frac{x^{2a+1}y^{2b+1}-((x-y)^{d-2}-xyg(x,y))^2}{(x-y)^{d-2}}=0,$$
where $d \geq 4$,  $g(x,y)=y^{d-3}h(x/y)$ and 
$$h(t)= \sum_{k=0,d-3}\frac{a_k}{k!}(t-1)^k,$$
with $a_0=1$, $a_1=a-\frac{1}{2}$ and $a_k=a_1(a_1-1) \cdots (a_1-k+1)$ for $k>1$, see \cite{FZ}.
Then it follows from \cite{Art}, Corollary 1, that $\pi_1(U)$ is abelian if $2a+1$ and $2b+1$ are relatively prime. It follows by Corollary \ref{corPi} that all these curves are either free or nearly free.
In \cite{DStFD}, Example 4.3 we have checked by a direct computation that these curves (without any condition on $2a+1$ and $2b+1$)  are actually free for $5 \leq d \leq 10$. The curve obtained for $d=4$ is almost free.
\end{ex}

To get futher examples, including curves of odd degree, we recall now the classification of unicuspidal rational curves with a unique Puiseux pair, see Theorem 1.1 in \cite{FLMN}.

\begin{thm}
\label{thmClass}

Let $a_i$ be the Fibonacci numbers with $a_0=0$, $a_1=1$, $a_{j+2}=a_{j+1}+a_{j}$.
A Puiseux pair $(a,b)$ can be realized by a unicuspidal rational curve of degree $d\geq 3$ if and only if the triple $(a,b,d)$ occurs in the following list.
\begin{enumerate}

\item $(d-1,d,d)$;

\item $(d/2,2d-1,d)$ with $d$ even;

\item $(a_{j-2}^2, a_{j}^2, a_{j-2} a_{j})$ with $j \geq 5$ odd;

\item $(a_{j-2}, a_{j+2}, a_{j})$ with $j \geq 5$ odd;

\item $(3,22,8)$;

\item $(6,43,16)$.

\end{enumerate}
\end{thm}

\begin{cor}
\label{corTop}
A unicuspidal rational curve with a unique Puiseux pair not of the type $(3)$ above with $d=a_{j-2} a_{j}$ odd is either free or nearly free.
\end{cor}

Numerical experiments suggest that the unicuspidal rational curves with a unique Puiseux pair  of the type (3) above with $d=a_{j-2} a_{j}$ odd is also either free or nearly free, but our method of proof does not work in this case.

\proof

Recall that a cusp with a unique Puiseux pair $(a,b)$ has the same monodromy eigenvalues as the cusp $u^a+v^b=0$. It follows that the order $k_x$ of any such eigenvalue $\lambda_x$ should be a divisor of $ab$, but not a divisor of either $a$ or $b$. This remark combined with 
 Theorem \ref{thm2} settles the cases (1), (2), (5) and (6) in Theorem \ref{thmClass}, as well as the case (3) when the degree $d$ is even.

For the case (4), we use the Catalan's identity
$$ a_j^2-a_{j-2}a_{j+2}=(-1)^j$$
and conclude that $k_x=1$ for any  eigenvalue $\lambda_x$ which is a $d$-th root of unity.

\endproof

The same idea as that in the proof of  Theorem \ref{thm2} gives the following, to be compared with Corollary \ref{corC1}.

\begin{prop} \label{bound} 
Let $C:f=0$ be a rational cuspidal curve of degree $d$. Then $N(f)_k=0$ for $k \leq d-3$ or $k\geq 2d-3$ and $st(f) \leq 2d-3$.

\end{prop} 

Propositions \ref{ST1} and \ref{ST2} below show that this vanishing result is sharp.

\proof
As in the proof of  Theorem \ref{thm2} we have an injection
$$N(f)_{2d-3} \to H^2(F, \C)_{1},$$
and $H^2(F, \C)_{1}=H^2(U, \C)$, where $U= \PP^2 \setminus C$. Since $E(U)=E(\PP^2)-E(C)=3-2=1$ and $b_0(U)=1$, $b_1(U)=0$ (since $C$ is irreducible), it follows that $b_2(U)=0$.
Hence $N(f)_{2d-3}=0$ and using the Lefschetz type property of the Jacobian module $N(f)$, see 
Corollary 4.3 in \cite{DPop}, it follows that $N(f)_{k}=0$ for any $k \geq 2d-3$.
To end the proof of the vanishing claim, it is enough to use the self duality of the graded module $N(f)$, see \cite{DS1}
or \cite{Se}, which yields $\dim N(f)_k=\dim N(f)_{T-k}$.

To prove the claim on the stability threshold $st(f)$, one uses the formula \eqref{n(f)} and Corollary 8 in \cite{CD}.

\endproof

\section{Further examples of nearly free and almost free divisors}

Consider first the simplest rational cuspidal curve of degree $d$ (it is the unique  rational curve up to projective equivalence with a unique cusp which is weighted homogeneous), namely
$C_d:f_d=x^d+y^{d-1}z=0$. 

\begin{prop} \label{ST1} 
The curve $C_d:f_d=x^d+y^{d-1}z=0$ is almost free for $d=2$ (when it is a smooth conic) and it is nearly free for $d \geq 3$, when it is a rational unicuspidal curve with $\mu(C)=\tau(C)=(d-1)(d-2).$ Moreover $(d_1,d_2)=(1,d-1)$ and $b=1$. In addition $\dim N(f)_k=1$ for $d-2 \leq k \leq 2d-4$ and $N(f)_k=0$ otherwise.
\end{prop} 

\proof

The polynomial $f_d$ is of Sebatiani-Thom type, i.e. the variables a separated in two groups. It follows that $M(f)=M(x^d) \otimes M(y^{d-1}z)$. 
 Moreover, we have $N(f)= N(x^d)\otimes
N(y^{d-1}z)$, and hence the graded module $N(f)$ has the following monomial basis $y^{d-2},xy^{d-2},..., x^{d-2}y^{d-2}$. This completes the proof.

\endproof

Exactly the same approach gives a similar result for a family having two singularities.

\begin{prop} \label{ST2} The curve
$C_d:f_d=x^d-y^{d-k}z^k=0$ for $d\geq 4$ and $2 \leq k \leq d/2$  is nearly free, has two singularities and $\tau(C)=(d-1)(d-2)$.
In addition, $(d_1,d_2)=(1,d-1)$, $b=1$,  $\dim N(f)_j=1$ for $d-2 \leq j \leq 2d-4$ and $N(f)_j=0$ otherwise. The number of irreducible components is given by the greatest common divisor $g.c.d. (d,k)$ and each of them is a rational curve.

\end{prop}

\begin{rk}
\label{rkdef}

\medskip

\noindent (i) Note that the family of irreducible free divisors
$$C_d: f_d=y^{d-1}z+x^d+ax^2y^{d-2}+bxy^{d-1}+cy^{d}=0, \  \  \ a \ne 0$$
described in \cite{ST} is a deformation of the curves in Proposition \ref{ST1}.
In the same spirit, the nearly free curve $f=x^7+y^5z^2$ with $\dim N(f)=6$ has two deformations which are almost free, namely $f'=f+x^3y^4$ and $f''=f+x^5y^2$.

\noindent (ii) The construction in Propositions \ref{ST1} and \ref{ST2} cannot be extended in general to curves given by $x^d+g(y,z)=0$, where $g(y,z)$ has at least 3 distinct linear factors. For example, the curve $f=x^7+y^2z^3(y+z)^2$ has $\dim N(f)_{[T/2]}=4$.

\end{rk}

\begin{ex} \label{ex2} The rational cuspidal curves of degree $d=6$  have been classified by Fenske, see \cite{F} and are either free or nearly free by Theorem \ref{thm2}. We discuss these curves in more detail now, and say which of them are free (resp. almost or nearly free), a fact which cannot be obtained from  Theorem \ref{thm2}, but usually by  using a computer algebra software.

When the curve has a unique cusp, there three classes of such curves.
The class $C_1$ contains the the curve $C:x^6-y^5z=0$ as the type (a), and 3 other additional families of curves denoted by (b), (c) and (d), which are deformations of the curve  $C$. The type (a) is nearly free with $(d_1,d_2,d_3)=(1,5,5)$ as seen in Proposition \ref{ST1}, and all the others are almost free with $(d_1,d_2,d_3)=(3,3,3)$.

The class $C_2$ contains  two types denoted by (a) and (b), which are suitable deformations of the curve $C:x^6-y^4z^2=0$. The families in both type (a)  and type (b) are almost free
with $(d_1,d_2,d_3)=(3,3,3)$. 

The class $C_3$ contains  three families denoted by (a), (b) and (c), which are suitable deformations of the curve $C:x^6-y^3z^3=0$. The family in the type (a) has $\tau=19$ and it is  free. The families in the types (b) and (c) are almost free with $(d_1,d_2,d_3)=(3,3,3)$.

When the curve has two cusps, there are six classes denoted by $C_4, C_5,C_6,C_7,C_8$ and $C_9$. Some of them, namely $C_4 (b)$, $C_5$, $C_6$ and $C_9$ are free divisors, see \cite{DStFD}, Example 4.2. The remaining three classes
$$C_4 (a): f=-a^3x^4y^2-2a^2x^3y^3-ax^2y^4-y^6+2a^2x^4yz+2ax^3y^2z+3xy^4z$$
$$-ax^4z^2-3x^2y^2z^2+
x^3z^3,$$
with $a \ne 0$, 
$$C_7: f=x^3y^3-A^2x^2y^4-2Axy^5-y^6+2Ax^3y^2z+2x^2y^3z-x^4z^2 $$
with $A \ne 0$ and 
$$C_8:  f=(1-2A)xy^5-y^6+2Ax^3y^2z+2x^2y^3z-x^4z^2-A^2x^2y^4$$
with $A \ne 0$ and  $A \ne 1/2$ are all almost free, with $(d_1,d_2,d_3)=(3,3,3)$.

\end{ex}

\begin{ex} \label{prop0}

 The rational cuspidal curves of degree $d$ and type $(d,d-2)$ with a unique cusp have been classified by Sakai and Tono, see \cite{SaTo}. 
Their result is recalled as Proposition 3.1 in  \cite{DStFD}. 
These curves  depend on parameters, have even degrees and hence are either free or nearly free by Theorem \ref{thm2}. 

Taking the simpliest normal forms in this class we get the  series of curves
$$C_d: f_d= (y^kz+x^{k+1})^2-xy^{2k+1}$$
of even degree $d=2k+2$.

These curves 
$C_d: f_d=0$
are almost free  for all $k \geq 2$,  $(d_1,d_2)=( k+1,k+1)$, $b=-k+1$, $ct=st= 3k+1$, $\tau=3k(k+1)$.

As explained in Remark \ref{rkreso}, in order to construct the corresponding resolution we need the following.

\noindent (i) the generating relations for $AR(f_d)$, which in this case are
$$r_1=(2(2k+1)y^kz, 2x^ky,(2k+1)y^{k+1}-2mx^kz),$$
$$r_2=(-2k(2k+1)xy^{k-1}z,2(k+1)(2k+1)x^{k+1}+2my^kz, (k+1)(2k+1)^2xy^k-2kmy^{k-1}z^2),$$
$$r_3=(-(2k+1)xy^k,y^{k+1}, (k+1)(2k+1)x^{k+1}-ky^kz),$$
where $m=2k^2+4k+1$.

\noindent (ii) the relation among $r_1$, $r_2$ and $r_3$, which in this case is
$$R=(-(k+1)(2k+1)x, y, -2mz).$$
\end{ex}

\begin{ex} \label{ex3} The rational cuspidal curves of degree $d$ and type $(d,d-2)$ with two cusps have been classified by Sakai and Tono, see \cite{SaTo}, who found three classes. 
Their result is recalled as Proposition 3.2 in  \cite{DStFD}, the classes corresponding to the cases (i), (ii) and (iii) in the statement. 

\medskip

The curves of class (i) for $5 \leq d \leq 15$ (and conjecturally for any $d\geq 5$ are free divisors, see \cite{DStFD}, Example 4.1. For $d$ even, these curves  are either free or nearly free by Theorem \ref{thm2}.

\medskip

The curves in class (ii) depend on parameters, the definig equation being 
$$C_d:f_d=(y^{k-1}z+\sum_{i=2,k}a_ix^i y^{k-i})^2y - x^{2k+1}=0,$$
where $d=2k+1 \geq 5$.
For one choice of these parameters, namely $a_k=1$ and the other $a_i=0$, one gets an infinite series of free divisors, see \cite{DStFD}, Theorem 4.6. The choices
$f_d=(y^{k-1}z+x^2y^{k-2})^2y-x^{2k+1}$, 
$f_d=(y^{k-1}z+x^{k-1}y)^2y-x^{2k+1}$, and 
$f_d=(y^{k-1}z+x^{k-1}y+x^k)^2y-x^{2k+1}$ also give free divisors for $k \geq 3$.
Note that the choice of parameters $a_i=0$ for all $i$ corresponds to Example \ref{ST2} above, case $d$ odd, hence we get a nearly free divisor.

\medskip

The curves in case (iii) also depend on parameters, have even degree and hence are either free or nearly free by Theorem \ref{thm2}. Taking the simpliest normal forms in this class we get the double series of curves
\begin{equation} \label{cjk}
C_{j,k}: f=(y^{k+j}z+x^{k+j+1})^2-x^{2j+1}y^{2k+1}=0.
\end{equation} 
Note that the curve $C_{j,k}$ has a cusp $A_{2j}$ for $1 \leq j\leq (d-2)/2$.
These curves are almost free for $2\leq k+j \leq 10$ (and conjecturally for all $k+j \geq 2$),
with $d=2(k+j)+2$, $d_1=d_2=d_3= k+j+1$, $b=-(d-4)/2$, $ct=st= (3d-4)/2$, $\tau=3d(d-2)/4$. In view of  Corollary \ref{corA}, it is enough to find the corresponding relations $r_1,r_2,r_3$ and $R$.

\end{ex}

\end{document}